\documentclass[12pt]{article}

\usepackage[utf8]{inputenc}

\usepackage{amsfonts}
\usepackage{amssymb}
\usepackage{amsmath}
\usepackage{amsthm}
\usepackage{amscd}
\usepackage{mathrsfs}
\usepackage[T2A]{fontenc}

\usepackage[russian]{babel}

\usepackage{amsthm}

\def \le {\leqslant}
\def \ge {\geqslant}

\theoremstyle{plain}

\newtheorem{lem}{Lemma}[section]
\newtheorem{foll}{Corollary}[section]
\newtheorem{theor}{Theorem}

\newenvironment{solve}{\begin{proof}[Proof]}{\end{proof}}
\addto\captionsrussian{}
\addto\captionsrussian{}

\newcommand\blfootnote[1]{%
  \begingroup
  \renewcommand\thefootnote{}\footnote{#1}%
  \addtocounter{footnote}{-1}%
  \endgroup
}

\begin{document}
\blfootnote{
                    Research is supported by RNF grant No. 14-11-00433\\
                    The author is a Young Russian Mathematics award winner and would like to thank its sponsors and jury.}
\centerline{\Large \bf Admissible endpoints of gaps in the Lagrange spectrum}
\vspace*{3mm}
\centerline{\Large Dmitry Gayfulin\footnote{ Research is supported by RNF grant No. 14-11-00433}}
\vspace*{3mm}

\vspace*{8mm}
\begin{abstract}
For any real number $\alpha$ define the Lagrange constant $\mu(\alpha)$ by
$$
\mu^{-1}(\alpha)=\liminf_{p\in\mathbb{Z}, q\in\mathbb{N}} |q(q\alpha-p)|.
$$

The set of all values taken by $\mu(\alpha)$ as $\alpha$ varies is called the \textit{Lagrange spectrum} $\mathbb{L}$. Irrational $\alpha$ is called attainable if the inequality 
$$
\biggl|\alpha -\frac{p}{q}\biggr|\le\frac{1}{\mu(\alpha)q^2}
$$
holds for infinitely many integers $p$ and $q$. Throughout this paper we will call a real number $\lambda\in\mathbb{L}$ \textit{admissible} if there exists an irrational attainable $\alpha$ such that $\mu(\alpha)=\lambda$. In our previous paper \cite{My} we constructed an example of not admissible element in the Lagrange spectrum. In the present paper, we give a necessary and sufficient condition of admissibility of a Lagrange spectrum element. We also prove that all elements of the sequence of $\mathbb{L}$ left endpoints $\alpha_n^*, n\ge 2$, first considered by Gbur, are not admissible.
\end{abstract}

\section{Introduction}
The Lagrange spectrum $\mathbb{L}$ is usually defined as the set of all values of the Lagrange constants 
$$
\mu(\alpha)=\biggl(\liminf_{p\in\mathbb{Z}, q\in\mathbb{N}} |q(q\alpha-p)|\biggr)^{-1}
$$
as $\alpha$ runs through the set of irrational numbers. Consider the continued fraction expansion of $\alpha$
$$
\alpha=[a_0;a_1,a_2,\ldots,a_n,\ldots].
$$
For any positive integer $i$ define 
$$
\lambda_i(\alpha)=[a_i;a_{i+1},a_{i+2},\ldots]+[0;a_{i-1},a_{i-2},\ldots,a_1],
$$
It is well known fact that
\begin{equation}
\label{lalimsup}
\limsup\lambda_i(\alpha)=\mu(\alpha).
\end{equation}
The equation (\ref{lalimsup}) provides an equivalent definition of the Lagrange constant $\mu(\alpha)$.

The following properties of $\mathbb{L}$ are well known. The Lagrange spectrum is a closed set \cite{Cusick1975} with minimal point $\sqrt{5}$. All the numbers of $\mathbb{L}$ which  are less than $3$ form a discrete set. It is well known fact that the Lagrange spectrum contains all elements over $\sqrt{21}$ (see\cite{Freiman1},\cite{Schecker}). The complement of $\mathbb{L}$ is a countable union of \textit{maximal gaps} of the spectrum. The maximal gaps are open intervals $(a,b)$ such that $(a,b)\cap\mathbb{L}=\emptyset$, but $a$ and $b$ both lie in the Lagrange spectrum. There are infinitely many gaps in the non-discrete part of the Lagrange spectrum \cite{Gbur}.

 Let $\alpha$ be an arbitrary irrational number. If the inequality
\begin{equation}
\label{attdef}
\biggl|\alpha -\frac{p}{q}\biggr|\le\frac{1}{\mu(\alpha) q^2}
\end{equation}
has infinitely many solutions for integer $p$ and $q$, we call $\alpha$ \textit{attainable}. This definition was first given by Malyshev in \cite{Malyshev}. One can easily see\cite{My} that $\alpha$ is attainable if and only if $\lambda_i(\alpha)\ge\mu(\alpha)$ for infinitely many indices $i$. We also call a real number $\lambda\in\mathbb{L}$ \textit{admissible} if there exists an irrational attainable number $\alpha$ such that $\mu(\alpha)=\lambda$.

Let $B$ denote a doubly infinite sequence of positive integers 
$$
B=(\ldots,b_{-n},\ldots b_{-1}, b_0, b_1,\ldots,b_n,\ldots).
$$
For an arbitrary integer $i$ define
$$
\lambda_i(B)=[b_i;b_{i-1},\ldots]+[0;b_{i+1},b_{i+2},\ldots].
$$
We will call a doubly infinite sequence $B$  purely periodic if there exists a finite sequence $P$ such that $B=(\overline{P})$.
A doubly infinite sequence $B$ is called eventually periodic if there exist $3$ finite sequences $P_l,R,P_r$ such that $B=(\overline{P_l},R,\overline{P_r})$.
One can also consider an equivalent definition of the Lagrange spectrum using the doubly infinite sequences. Denote
\begin{equation}
\label{LMdef}
L(B)=\limsup\limits_{i\to\infty}\lambda_i(B),\quad M(B)=\sup\lambda_i(B). 
\end{equation}
The Lagrange spectrum $\mathbb{L}$ is exactly the set 
of values taken by $L(B)$ as $B$ runs through the set of doubly infinite sequences of positive integers. The set of values taken by $M(B)$ is called the Markoff spectrum. We will denote this set by $\mathbb{M}$.

We will call a doubly infinite sequence $B$ \textit{weakly associated} with an irrational number $\alpha=[a_0;a_1,\ldots,a_n,\ldots]$ if the following condition holds:
\begin{enumerate}
\item{For any natural $i$ the pattern $(b_{-i},b_{-i+1},\ldots,b_0,\ldots,b_i)$ occurs in the sequence $a_1,a_2,\ldots,a_n,\ldots$ infinitely many times.}

We will call $B$ \textit{strongly associated} with $\alpha$ if, additionally,
\item{$\mu(\alpha)=\lambda_0(B)=M(B)$.}
\end{enumerate}
One can easily see that if $B$ is weakly associated with $\alpha$ then $\mu(\alpha)\ge M(B)$. As we will show in Section \ref{secevper}, if $\alpha$ has limited partial quotients, it has at least one strongly associated sequence.

\section{Results of paper  \cite{My}}

\textbf{Theorem I}\\
\textit{The quadratic irrationality $\lambda_0=[3;3,3,2,1,\overline{1,2}]+[0;2,1,\overline{1,2}]$ belongs to $\mathbb{L}$, but if $\alpha$ is such that $\mu(\alpha)=\lambda_0$ then $\alpha$ is not attainable.}\\
\textbf{Theorem II}\\
\textit{If $\lambda\in\mathbb{L}$ is not a left endpoint of some maximal gap in the Lagrange spectrum then there exists an attainable $\alpha$ such that $\mu(\alpha)=\lambda$.}\\
One can easily formulate these theorems using the concept of admissible numbers, introduced behind.
\\
\textbf{Theorem I'}\\
\textit{The quadratic irrationality $\lambda_0=[3;3,3,2,1,\overline{1,2}]+[0;2,1,\overline{1,2}]$ belongs to $\mathbb{L}$, but is not admissible.}\\
\textbf{Theorem II'}\\
\textit{If $\lambda\in\mathbb{L}$ is not a left endpoint of some maximal gap in the Lagrange spectrum then $\lambda$ is an admissible number.}\\

\section{Main results}
Our first theorem is a small generalization of Theorem 3 in \cite{My}. The proof will be quite similar and use some lemmas from \cite{My}.
\begin{theor}
\label{evper}
Let $a$ be a left endpoint of a gap $(a,b)$ in the Lagrange spectrum and $\alpha$ be an irrational number such that $\mu(\alpha)=a$. Consider a doubly infinite sequence $B$ strongly associated with $\alpha$. Then $B$ is an eventually periodic (i.e. periodic at both sides) sequence.
\end{theor}
It follows from theorems I and II that there exist non-admissible elements in the Lagrange spectrum but all such numbers are left endpoints of some maximal gaps in $\mathbb{L}$. The following theorem gives a necessary and sufficient condition of admissibility of a Lagrange spectrum element.
\begin{theor}
\label{main}
A Lagrange spectrum left endpoint $a$ is admissible if and only if there exists a quadratic irrationality $\alpha$ such that $\mu(\alpha)=a$. 
\end{theor}

Of course, every quadratic irrationality is strongly associated with the unique sequence, which is purely periodic. Therefore Theorem \ref{main} is equivalent to the following statement.
\begin{foll}
\label{mainfoll2}
A Lagrange spectrum left endpoint $a$ is not admissible if and only if there does not exist a purely periodic sequence $B$ such that  $\lambda_0(B)=M(B)=a$.
\end{foll}
Theorem \ref{main} provides an instrument to verify not admissible points in $\mathbb{L}$. Define 
$$
\alpha_n^*=2+[0;\underbrace{1,\ldots,1}_{2n-2},\overline{2,2,1,2}]+[0;\underbrace{1,\ldots,1}_{2n-1},2,\underbrace{1,\ldots,1}_{2n-2},\overline{2,2,1,2}]
$$ and 
$$
\beta_n=2+2[0;\overline{\underbrace{1,\ldots,1}_{2n},2}].
$$ The fact that $(\alpha_n^*,\beta_n)$ is the maximal gap in the Markoff spectrum was proved by Gbur in \cite{Gbur}. It is easy to show that $\alpha_n^*$ and $\beta_n$ belong to $\mathbb{L}$. We will do it in the section \ref{Gburproof}. Hence, as $\mathbb{L}\subset\mathbb{M}$ \cite{Cusick1975}, the interval $(\alpha_n^*,\beta_n)$ is the maximal gap in $\mathbb{L}$ too.
\begin{theor}
\label{GburSer}
For any integer $n\ge 2$ the irrational number $\alpha_n^*$ is not admissible.
\end{theor}
One can easily see that $\alpha_1^*=2+[0;\overline{2,2,1,2}]+[0;1,2,\overline{2,2,1,2}]=\mu([0;\overline{2,2,1,2}])=M(\overline{2,2,1,2})$. Thus, $\alpha_1^*$ is an admissible number by Theorem \ref{main}.

\section{Proof of Theorem \ref{evper}}
\label{secevper}
The following statement is well known. See proof in (\cite{Cusick}, Ch. 1, Lemma 6).
\begin{lem}
Let $A=\ldots,a_{-1},a_0,a_1,\ldots$ be any doubly infinite sequence. If M(A) is finite, then there exists a doubly infinite sequence $B$ such that $M(A)=M(B)=\lambda_0(B)$.
\end{lem}
Using the same argument for the sequence $A=(a_1,a_2,\ldots,a_n,\ldots)$, one can easily show that
\begin{lem}
Let $\alpha=[0;a_1,\ldots,a_n,\ldots]$ be an arbitrary irrational number and $a_i<c \ \forall i\in\mathbb{N}$ for some positive real number $c$. Then there exists a doubly infinite sequence $B$ which is strongly associated with $\alpha$.
\end{lem}
As  $\alpha\le\sqrt{21}$, all elements of $B$ are bounded by $4$. For any natural $n$ denote $\varepsilon_n=2^{-(n-1)}, \delta_n=5^{-2(n+2)}$. We need the following lemmas from \cite{My}.
\begin{lem}
\label{comp}
Suppose $\alpha=[a_0;a_1,\ldots,a_n,b_1,\ldots]$ and  $\beta=[a_0;a_1,\ldots,a_n,c_1,\ldots]$, where $n\ge 0$, $a_0$ is an integer, $a_1,\ldots,a_n,b_1,b_2,\ldots,c_1,c_2,\ldots$ are positive integers bounded by $4$ with $b_1\ne c_1$. Then for $n$ odd, $\alpha>\beta$ if and only if $b_1>c_1$; for $n$ even, $\alpha>\beta$ if and only if $b_1<c_1$. Also, 
$$
\delta_n<|\alpha-\beta|<\varepsilon_n.
$$
\end{lem}
\begin{lem}
\label{smallshift}
Let $\gamma=[0;c_1,c_2,\ldots,c_N,\ldots]$ and $\gamma'=[0;c'_1,c'_2,\ldots,c'_N,\ldots]$ be two irrational numbers with partial quotients not exceeding $4$. Suppose that every sequence of partial quotients of length $2n+1$ which occurs in the sequence $(c'_1,c'_2,\ldots,c'_N,\ldots)$ infinitely many times also occurs in the sequence $(c_1,c_2,\ldots,c_N,\ldots)$ infinitely many times. Then $\mu(\gamma')<\mu(\gamma)+2\varepsilon_{n}$.
\end{lem}
The following technical lemma was formulated in \cite{My} for $N=(2n+1)(4^{2n+1}+1)$ and the proof was incorrect. However, this is not crucial for the results of the paper \cite{My} as we just need $N$ to be bounded from above by some growing function of $n$. In this paper, we give a new version of lemma with correct proof.
\begin{lem}
\label{Dirichlet}
Let $n$ be an arbitrary positive integer. Denote $N=N(n)=(2n+2)(4^{2n+2}+1)$. If $b_1, b_2,\ldots,b_N$ is an arbitrary integer sequence of length $N$ such that $1\le a_i\le 4$ for all $1\le i\le N$, then there exist two integers $n_1, n_2$ such that $b_{n_1+i}=b_{n_2+i}$ for all $0\le i\le 2n+1$ and $n_1\equiv n_2 \pmod 2$.
\end{lem}
\begin{solve}
There exist only $4^{2n+2}$ distinct sequences of length $2n+2$ with elements $1,2,3,4$. Consider $4^{2n+2}+1$ sequences: $(a_1, \ldots,a_{2n+2}), (a_{2n+3}, \ldots,a_{4n+4}),\ldots, (a_{(2n+2)4^{2n+2}+1},\ldots, a_{(2n+2)4^{2n+2}+2n+2})$. Dirchlet's principle implies that there exist two coinciding sequences among them. Denote these sequences by $(a_{n_1},\ldots,a_{n_1+2n+1})$ and $(a_{n_2},\ldots,a_{n_2+2n+1})$. Note that the index of the first element of each sequence is odd, hence $n_1\equiv n_2 \equiv 1\pmod 2$, that finishes the proof. 
\end{solve}
If $n_1\equiv n_2 \pmod 2$ then the sequence $(a_{n_1},a_{n_1+1},\ldots,a_{n_2-1})$ has even length. This fact will be useful in our argument.
\begin{lem}
\label{forSurgery}
Let $B$ be an arbitrary integer sequence of even length. Let $A$ be an arbitrary finite integer sequence and $C$ --- an arbitrary non-periodic infinite sequence. Then 
\begin{equation}
\label{leabbc}
\min([0;A,B,B,C], [0;A,C])<[0;A,B,C]<\max([0;A,B,B,C], [0;A,C])
\end{equation}

\begin{solve}
As the sequence $C$ is non-periodic, the continued fractions in (\ref{leabbc}) are not equal. Without loss of generality, one can say that the sequence $A$ is empty. Suppose that
$$
[0;B,C]>[0;B,B,C].
$$
As the length of $B$ is even, one can see that $[0;C]>[0;B,C]$, which is exactly the right part of the inequality (\ref{leabbc}). The case when $[0;B,C]<[0;B,B,C]$ is treated in exactly the same way.
\end{solve}

\end{lem}
\begin{lem}
\label{Surgery}
Let $\gamma=[0;b_1,b_2,\ldots,b_N,\ldots]$ be an arbitrary irrational number, not a quadratic irrationality. Consider the sequence $B_N=(b_1,b_2,\ldots,b_N)$ and define two numbers $n_1$ and $n_2$ from Lemma \ref{Dirichlet}. Define two new sequences of positive integers
\begin{equation*}
\begin{split}
B_N^1=(b_1,b_2,\ldots, b_{n_1-1},b_{n_2},b_{n_2+1},\ldots,b_N),\\
B_N^2=(b_1,b_2,\ldots, b_{n_1-1},b_{n_1},\ldots,b_{n_2-1},b_{n_1},\ldots,b_{n_2-1},b_{n_2},b_{n_2+1},\ldots,b_N).
\end{split}
\end{equation*}
Let us also define two new irrational numbers:
\begin{equation*}
\begin{split}
\gamma^1=[0;b_1,b_2,\ldots, b_{n_1-1},b_{n_2},b_{n_2+1},\ldots,b_N,b_{N+1}\ldots]=[0;B_N^1,b_{N+1},\ldots],\\
\gamma^2=[0;b_1,b_2,\ldots, b_{n_1-1},b_{n_1},\ldots,b_{n_2-1},b_{n_1},\ldots,b_{n_2-1},b_{n_2},b_{n_2+1},\ldots,b_N,\ldots]=[0;B_N^2,b_{N+1},\ldots].
\end{split}
\end{equation*}
Then $\max(\gamma^1, \gamma^2)>\gamma$.
\end{lem}
\begin{solve}
We apply Lemma \ref{forSurgery} for $A=(b_1,b_2,\ldots,b_{n_1-1}), B=(b_{n_1},b_{n_1+1},\ldots,b_{n_2-1}), C=(b_{n_2},b_{n_2+1},\ldots)$. Here $\gamma=[0;A,B,C], \gamma^1=[0;A,C], \gamma^2=[0;A,B,B,C]$. Note that as $\gamma$ is not a not a quadratic irrationality, the sequence $C$ is not periodic.
\end{solve}
Now we are ready to prove Theorem \ref{evper}.
\begin{solve}
Suppose that $B$ is not periodic at right side. Consider a growing sequence of indices $k(j)$ such that for any natural $j$ the sequence $(a_{k(j)-j},\ldots,a_{k(j)},\ldots,a_{k(j)+j})$ coincides with the sequence $(b_{-j},\ldots,b_0,\ldots,b_j)$. Of course,
$$
\lim\limits_{j\to\infty}\lambda_{k(j)}(\alpha)=\lambda_0(B)=\mu(\alpha).
$$
Without loss of generality, on can say that $k(j+1)-k(j)\to\infty$ as $j\to\infty$. Consider an even $n$ such that $\varepsilon_n<\frac{b-a}{2}$ and $N=N(n)$ as defined in Lemma \ref{Dirichlet}. Define $n_1<n_2$ from  Lemma \ref{Dirichlet} for the sequence $(b_1,\ldots,b_N)$. As $B$ is not periodic to the right, define a minimal positive integer $r$ such that $b_{n_1+r}\ne b_{n_2+r}$. Consider the sequences $B_N^1, B_N^2$ and the continued fractions $\gamma_1, \gamma_2$ from Lemma \ref{Surgery} applied to the continued fraction $[0;b_1,\ldots,b_n\ldots]=\gamma$. If $\gamma_2>\gamma$, define $g=2$, otherwise we put $g=1$. Consider the doubly infinite sequence $B'=(\ldots,b_{-n},b_0,B_N^g,b_{N+1},\ldots)$. Note that
$$
a=\lambda_0(B)<\lambda_0(B')<a+\varepsilon_n<b.
$$
Consider the corresponding continued fraction $\alpha'$ which is obtained from the continued fraction $\alpha$ by replacing every segment $(a_{k(j)},\ldots,a_{k(j)+N})=(a_{k(j)},B_N)$ by the segment $(a_{k(j)},B_N^g)$ for every $j\ge n_2+r$. One can easily see that $\alpha'$ and $\alpha$ satisfy the condition of lemma \ref{smallshift} and hence $\mu(\alpha')<\mu(\alpha)+2\varepsilon_{n}$. But as $\mu(\alpha)+2\varepsilon_{n}<b$ and $(a,b)$ is the gap in $\mathbb{L}$, we have 
\begin{equation}
\label{lessthana}
\mu(\alpha')\le\mu(\alpha)=a.
\end{equation}
From the other hand, one can easily see that the sequence $B'$ is weakly associated with $\alpha'$. This means that
$$
\mu(\alpha')\ge M(B)\ge\lambda_0(B')>\lambda_0(B)=a.
$$
We obtain a contradiction with (\ref{lessthana}). The case when $B$ is not periodic at left side is considered in exactly the same way. The theorem is proved.
\end{solve}
\section{Proof of Theorem \ref{main}}
The following Lemma from \cite{My} immediately implies the $\Leftarrow$ part of the statement of Theorem $\ref{main}$.
\begin{lem}
Consider an arbitrary point $a$ in the Lagrange spectrum. If there exists a quadratic irrationality $\gamma$ such that $\mu(\gamma)=a$, then $a$ is admissible.
\end{lem}
Now it is sufficient to prove that if $a$ is an admissible left endpoint of the Lagrange spectrum, then there exists a quadratic irrationality $\alpha$ such that  $\mu(\alpha)=a$.
\begin{solve}
Let $a$ be an admissible left endpoint of the Lagrange spectrum. Let $\alpha=[a_0;a_1,\ldots,a_n,\ldots]$ be an irrational number such that $\mu(\alpha)=a$. Suppose that $\alpha$ is attainable, but not a quadratic irrationality. Let $k(j)$ be a growing sequence of indices such that
\begin{equation}
\label{isatt}
\lambda_{k(j)}(\alpha)\ge\mu(\alpha).
\end{equation}
Of course,
$$
\lim\limits_{j\to\infty}\lambda_{k(j)}(\alpha)=\mu(\alpha).
$$
Consider a strongly associated with $\alpha$ sequence $B=(\ldots,b_{-n},\ldots b_{-1}, b_0, b_1,\ldots,b_n,\ldots)$ having the following property: the sequence $(b_{-i},\ldots,b_0,\ldots,b_i)$ coincides with the sequence $a_{k(j)-i},\ldots,a_{k(j)},\ldots,a_{k(j)+i}$ for infinitely many $j$-s. Theorem \ref{evper} implies that $B$ is eventually periodic. That is, there exist a positive integer $m$ and two finite sequences $L$ and $R$ such that
$$
B=(\overline{L},b_{-m},\ldots,b_0,\ldots,b_m,\overline{R}).
$$
 It follows from (\ref{isatt})  that one of the inequalities
\begin{equation}
\begin{split}
\label{gorightleft}
[a_{k(j)};a_{k(j+1)},\ldots]\ge[b_0;b_1,\ldots,b_m,\overline{R}], \text{or}\\
[0;a_{k(j-1)},\ldots,a_1]\ge[0;b_{-1},\ldots,b_{-m},\overline{L}]
\end{split}
\end{equation}
holds for infinitely many $j$-s. Note that $[a_{k(j)};a_{k(j+1)},\ldots]\ne[b_0;b_1,\ldots,b_m,\overline{R}]$, as $\alpha$ is not a quadratic irrationality and, of course, $[0;a_{k(j-1)},\ldots,a_1]\ne[0;b_{-1},\ldots,b_{-m},\overline{L}]$.
Suppose that
\begin{equation}
\label{goright}
[a_{k(j)};a_{k(j+1)},\ldots]>[b_0;b_1,\ldots,b_m,\overline{R}]
\end{equation}
for infinitely many $j$-s. Denote by $p$ the length of period $R$. Denote by $r(j)$ the minimal positive number such that $a_{k(j)+r(j)}\ne b_{r(j)}$. Without loss of generality, one can say that:
\begin{enumerate}
\item{$k(j+1)-k(j)-r(j)\to\infty$ as $j\to\infty$.}
\item{$[a_{k(j)};a_{k(j+1)},\ldots]>[b_0;b_1,\ldots,b_m,\overline{R}]$ for every $j\in\mathbb{N}$.}
\item{$[a_{k(j)};a_{k(j+1)},\ldots,a_{k(j)+m}]=[b_0;b_1,\ldots,b_m]$ for every $j\in\mathbb{N}$.}
\item{The sequence $(a_{k(j)-j},\ldots,a_{k(j)},\ldots,a_{k(j)+j})$ coincides with the sequence $(b_{-j},\ldots,b_0,\ldots,b_j)$ \\ for every $j \in\mathbb{N}$.}
\item{Period length $p$ is even.}
\end{enumerate}
 Denote by $t(j)$ the number of periods $P$ in the sequence $(b_{m+1},\ldots,b_{r(j)}).$ Of course, $t(j)=[\frac{r(j)-m}{p}]$ and $t(j)$ tends to infinity. Denote by $\alpha_n$ a continued fraction obtained from the continued fraction $\alpha=[a_0;a_1,\ldots,a_n,\ldots]$ as follows: if $t(j)>n$, then every pattern
$$
a_{k(j)},a_{k(j+1)},\ldots,a_{k(j)+m},\underbrace{R,\ldots,R}_{t(j)\ times},\ldots,a_{k(j)+r(j)}
$$
is replaced by the pattern
$$
a_{k(j)},a_{k(j+1)},\ldots,a_{k(j)+m},\underbrace{R,\ldots,R}_{n\ times},\ldots,a_{k(j)+r(j)}.
$$
Lemma \ref{comp} implies that since (\ref{goright}) 
$$
[a_{k(j)};a_{k(j+1)},\ldots,a_{k(j)+m},\underbrace{R,\ldots,R}_{t(j)\ times},\ldots,a_{k(j)+r(j)}]>[b_0;b_1,\ldots,b_m,\overline{R}]
$$
and the length of the period $R$ is even, one has
\begin{equation}
\label{mainineq}
[a_{k(j)};a_{k(j+1)},\ldots,a_{k(j)+m},\underbrace{R,\ldots,R}_{n\ times},\ldots,a_{k(j)+r(j)}]>[b_0;b_1,\ldots,b_m,\overline{R}]+\delta_{m+(n+1)p}.
\end{equation}
Since $k(j+1)-k(j)-r(j)\to\infty$ as $j\to\infty$ and the sequence $(a_{k(j)-j},\ldots,a_{k(j)})$ coincides with the sequence $(b_{-j},\ldots,b_0)\ \forall j \in\mathbb{N}$, one can easily see that $\mu(\alpha_n)\ge\mu(\alpha)+\delta_{m+(n+1)p}$.

Note that
\begin{equation}
\label{maineq}
\lim\limits_{n\to\infty}\mu(\alpha_n)=\mu(\alpha)=a.
\end{equation}
Indeed, every pattern of length $np$ which occurs  in the sequence of partial quotients of $\alpha$ infinitely many times, occurs in the sequence of partial quotients of $\alpha_n$ infinitely many times. Similarly, every pattern of length $np$ which occurs in the sequence of partial quotients of $\alpha_n$ infinitely many times, occurs in the sequence of partial quotients of $\alpha$ infinitely many times too. Then, by Lemma \ref{smallshift}, 
$$
|\mu(\alpha)-\mu(\alpha_n)|<2\varepsilon_{np}=2^{-np+2}\to 0
$$
as $n\to\infty$. We obtain a contradiction with the fact that $a$ is the left endpoint of the gap $(a,b)$ in the Lagrange spectrum. Indeed, the inequality (\ref{mainineq}) implies that $\mu(\alpha_n)>\mu(\alpha)\ \ \forall n\in\mathbb{N}$. In addition, the equality (\ref{maineq}) implies that there exists a positive integer $N$ such that for any $n>N$ one has $a=\mu(\alpha)<\mu(\alpha_n)<b$. 

If the inequality (\ref{goright}) does not hold infinitely many times, then the inequality
\begin{equation*}
\label{goleft}
[0;a_{k(j-1)},\ldots,a_1]>[0;b_{-1},\ldots,b_{-m},\overline{L}]
\end{equation*}
holds infinitely many times. This case is treated in exactly the same way. The theorem is proved.
\end{solve}
\section{Proof of Theorem \ref{GburSer}}
\label{Gburproof}
First of all, let us show that $(\alpha_n^*,\beta_n)$ is the maximal gap in $\mathbb{L}$. As 
$$
\beta_n=2+2[0;\overline{\underbrace{1,\ldots,1}_{2n},2}]=\mu([0;\overline{\underbrace{1,\ldots,1}_{2n},2}]),
$$
we have $\beta_n\in\mathbb{L}$. The proof of the fact that $\alpha_n^*\in\mathbb{L}$, when $n\ge 2$ is little more complicated. Recall that 
$$
\alpha_n^*=2+[0;\underbrace{1,\ldots,1}_{2n-2},\overline{2,2,1,2}]+[0;\underbrace{1,\ldots,1}_{2n-1},2,\underbrace{1,\ldots,1}_{2n-2},\overline{2,2,1,2}].
$$
Denote by $C_n(k)$ the following finite sequence of integers
\begin{equation}
\label{starmark}
C_n(k)=(\underbrace{2,1,2,2}_{k},\underbrace{1,\ldots,1}_{2n-2},2^{*},\underbrace{1,\ldots,1}_{2n-1},2^{**},\underbrace{1,\ldots,1}_{2n-2},\underbrace{2,2,1,2}_{k}).
\end{equation}
Denote by $\zeta_n$ the following infinite continued fraction:
$$
\zeta_n=[0;C_n(1),C_n(2),\ldots,C_n(k),\ldots].
$$
A little calculation shows that $\mu(\zeta_n)=\alpha_n^*$ and therefore $\alpha_n^*$ belongs to the Lagrange spectrum $\mathbb{L}$.
By (\cite{Gbur}, Lemma 4), $\alpha_n^*$ is a growing sequence. One can easily see that 
$$  
\lim\limits_{n\to\infty}\alpha_n^*=2+2[0;\overline{1}]=\sqrt{5}+1\approx 3.236.
$$
Thus, we have
\begin{equation}
\label{alphanineq}
\alpha_2^*\le\alpha_n^*<1+\sqrt{5}\quad \text{where}\ n\ge 2.
\end{equation}
The following lemma is compilation of lemmas 3 and 4 from \cite{Gbur}
\begin{lem}
\label{banned}
Consider a doubly infinite sequence $B=(\ldots,b_{-n},\ldots,b_{-1},b_0,b_1,\ldots,b_n,\ldots)$ such that $M(B)<\sqrt{5}+1$. Then all elements of $B$ are bounded by $2$ and $B$ does not contain patterns of the form $(2,1,2,1)$ and $(1,2,1,2)$.
\end{lem}
%\begin{solve}
%By lemma 3 in \cite{Gbur}We just note that $1+\sqrt{5}<2+[0;\overline{1,2,2,1,2,1,1,2}]+[0;\overline{1,1,2,1,2,2,1,2}]$. 
%\end{solve}
By Lemma \ref{secevper}, without loss of generality one can say that $M(B)=\lambda_0(B)$. Denote the continued fractions $[0;b_1,\ldots,b_n,\ldots]$ and $[0;b_{-1},\ldots,b_{-n},\ldots]$ by $x$ and $y$ respectively. Then
$$
M(B)=b_0+x+y.
$$
Without loss of generality one can say that $x\le y$. Now we need the following lemma from (\cite{Gbur}, Theorem 4(i)).
\begin{lem}
\label{firstelements}
Let $B$ be a doubly infinite sequence such that $M(B)=\lambda_0(B)$, then for all $n\ge 1$ we have
\begin{equation}
\label{maincompare}
\begin{split}
\beta_n\le M(B)=2+x+y\le \alpha_{n+1}^* \Leftrightarrow x=[0;\underbrace{1,\ldots,1}_{2n},2,\ldots]\ \text{and}\ y=[0;\underbrace{1,\ldots,1}_{2n},\ldots].\\
\end{split}
\end{equation}
\end{lem}
It also follows from (\cite{Gbur}, Theorem 4(ii)) that
\begin{equation}
\label{maincompare2}
 2+[0;\underbrace{1,\ldots,1}_{2n+1},\ldots]+[0;\underbrace{1,\ldots,1}_{2n+1},2,\ldots]<\sqrt{5}+1.
\end{equation}
Denote
\begin{equation}
\begin{split}
w_0=[0;\overline{2,1,2,2}],\ x_0=[0;\underbrace{1,\ldots,1}_{2n},\overline{2,2,1,2}]=[0;\underbrace{1,\ldots,1}_{2n},2+w_0],\\
y_0=[0;\underbrace{1,\ldots,1}_{2n+1},2,\underbrace{1,\ldots,1}_{2n},\overline{2,2,1,2}]=[0;\underbrace{1,\ldots,1}_{2n+1},2,\underbrace{1,\ldots,1}_{2n},2+w_0].
\end{split}
\end{equation}
\begin{lem}
\label{2212min}
Let $w=[0;a_1,a_2,\ldots,a_n,\ldots]$ be a continued fraction with elements equal to $1$ or $2$. Suppose that the sequence $(a_1,a_2,\ldots,a_n,\ldots)$ does not contain the pattern $(2,1,2,1)$. Then $w\ge w_0$
\end{lem}
\begin{solve}
Denote the elements of the continued fraction $w_0=[0;\overline{2,1,2,2}]$ by $[0;a'_1,\ldots,a'_m,\ldots]$. Denote by $r$ the minimal index such that $a_r\ne a'_r$. Suppose that $w<w_0$. Then either $r$ is odd, $a_r=2, a'_r=1$ or $r$ is even, $a_r=1, a'_r=2$. However $a'_r=2$ for any even $r$, thus the first case leads to a contradiction. Consider the second case. Of course, $r\ge 4$. Then $a'_{r-3}=a_{r-3}=2, a'_{r-2}=a_{r-2}=1, a'_{r-1}=a_{r-1}=2$. This means that $(a_{r-3},a_{r-2},a_{r-1},a_r)=(2,1,2,1)$ and we obtain a contradiction.
\end{solve}
\begin{lem}
\label{almostgbur}
If $B$ is strongly associated with $\alpha_{n+1}^*$, having 
$$
M(B)=\lambda_0(B)=2+x+y=\alpha_{n+1}^*.
$$ 
Then $x=x_0$ and $y=y_0$.
\end{lem}
\begin{solve}
By Lemma \ref{firstelements} $x=[0;\underbrace{1,\ldots,1}_{2n},2,\ldots]$ and  $y=[0;\underbrace{1,\ldots,1}_{2n},\ldots]$.
Note that $x\le x_0$. Indeed, 
$$
x=[0;\underbrace{1,\ldots,1}_{2n},2,b_{2n+2},\ldots]\le[0;\underbrace{1,\ldots,1}_{2n},2+w_0]=x_0 \Leftrightarrow [0;b_{2n+2},\ldots]\ge w_0.
$$
The last equality follows from lemmas \ref{2212min} and \ref{banned}.

Now suppose that 
$$
y=[0;b_{-1},\ldots,b_{-n},\ldots]>y_0.
$$
Denote the elements of the continued fraction $y_0$ by $[0;b'_1,\ldots,b'_m,\ldots]$. Denote by $r$ the minimal positive integer such that $b_{-r}\ne b'_r$. Lemma  \ref{firstelements} implies that $r>2n$. As $y>y_0,\ b_{-2n-1}\le b'_{2n+1}=1$, therefore $b_{-2n-1}=1$. Similarly, $b_{-2n-2}\ge b'_{2n+2}=2$,  Lemma \ref{banned} implies that $b_{-2n-2}=2$. Now we have
$$
y=[0;\underbrace{1,\ldots,1}_{2n+1},2,b_{-2n-3},\ldots]>[0;\underbrace{1,\ldots,1}_{2n+1},2,\underbrace{1,\ldots,1}_{2n},\overline{2,2,1,2}]=y_0. 
$$
Suppose that $r\le 4n+2$. Then $r$ is even, because otherwise $b_{-r}<b'_r=1$. Hence $2=b_{-r}>b'_r=1$. Note that $\lambda_{-2n-2}(B)\le\lambda_0(B)=M(B)$. As $r-2n-3<2n+1$, the inequality (\ref{maincompare2}) implies that
\begin{equation}
\label{111212}
\lambda_{-2n-2}(B)=2+[0;\underbrace{1,\ldots,1}_{2n+1},\ldots]+[0;\underbrace{1,\ldots,1}_{r-2n-3},2,\ldots]>1+\sqrt{5}.
\end{equation}
Thus, $r>4n+2$. Note that $r\ne 4n+3$, because otherwise  by (\ref{maincompare2}) one has
$$
\lambda_{-2n-2}(B)=2+[0;\underbrace{1,\ldots,1}_{2n+1},2,\ldots]+[0;\underbrace{1,\ldots,1}_{2n+1},\ldots]>1+\sqrt{5}.
$$
Hence, we have
$$
y=[0;\underbrace{1,\ldots,1}_{2n+1},2,\underbrace{1,\ldots,1}_{2n},2,b_{-4n-4},\ldots]>[0;\underbrace{1,\ldots,1}_{2n+1},2,\underbrace{1,\ldots,1}_{2n},2+w_0]=y_0  \Leftrightarrow [0;b_{-4n-4},\ldots] < w_0.
$$
We obtain a contradiction with Lemma  \ref{2212min}. The lemma is proved.
\end{solve}
In other words, 
$$ 
M(B)=\lambda_0(B)=\alpha_{n+1}^* \Leftrightarrow B=(\overline{2,1,2,2},\underbrace{1,\ldots,1}_{2n},2,\underbrace{1,\ldots,1}_{2n+1},2,\underbrace{1,\ldots,1}_{2n},\overline{2,2,1,2})
$$
Now the proof of Theorem \ref{GburSer} is quite simple.
\begin{solve}
Suppose that $\alpha^*_n$ is admissible for some $n\ge 2$. Consider an attainable number $\alpha$ such that $\mu({\alpha})=\alpha^*_n$. Theorem \ref{main} implies that $\alpha$ is a quadratic irrationality. Then $\alpha$ is associated with the unique purely periodic sequence $(\overline{P})$, where $P$ is the period of continued fraction expansion of $\alpha$. Lemma \ref{almostgbur} implies that 
\begin{equation}
\label{finalarg}
(\overline{P})=(\overline{2,1,2,2},\underbrace{1,\ldots,1}_{2n},2,\underbrace{1,\ldots,1}_{2n+1},2,\underbrace{1,\ldots,1}_{2n},\overline{2,2,1,2}),
\end{equation}
which is impossible because the sequence from the right part of (\ref{finalarg}) is not purely periodic. We obtain a contradiction and the theorem is proved.
\end{solve}
%Inequality (\ref{isatt}) implies that 
%$$
%[a_{k(j)};a_{k(j+1)}\ldots,a_{k(j)+r(j)}]>[b_0;b_1,\ldots]
%$$

Dmitry Gayfulin,\\
Steklov Mathematical Institute of Russian Academy of Sciences\\
ul. Gubkina, 8, Moscow, Russia, 119991\\
\textit{gayfulin@rambler.ru}
\end{document}